\documentclass[11pt]{smfart}

\usepackage{smfthm}
\usepackage{smfenum}
\usepackage{amssymb}
\usepackage{amscd}

\newcommand{\Qp}{\mathbf{Q}_p}
\newcommand{\Zp}{\mathbf{Z}_p}
\newcommand{\Cp}{\mathbf{C}_p}

\newcommand{\eps}{\varepsilon}
\newcommand{\ra}{\rightarrow}

\newcommand{\on}{\operatorname}
\newcommand{\OO}{\mathcal{O}}

\renewcommand{\hat}{\widehat}

\renewcommand{\phi}{\varphi}
\renewcommand{\projlim}{\varprojlim}
\newcommand{\ZZ}{\mathbf{Z}}
\newcommand{\NN}{\mathbf{N}}

\newcommand{\bnrig}[2]{\mathbf{B}^{\dagger #1}_{\mathrm{rig} #2}}

\newcommand{\btrigplus}[1]{\widetilde{\mathbf{B}}^{+}_{\mathrm{rig} #1}}

\newcommand{\dnrig}[1]{\mathbf{D}^{\dagger #1}_{\mathrm{rig}}}
\newcommand{\dhol}{\mathbf{D}^+_{\mathrm{rig}}}

\newcommand{\bmax}{\mathbf{B}_{\mathrm{max}}}

\newcommand{\bdR}{\mathbf{B}_{\mathrm{dR}}}

\newcommand{\bplus}{\mathbf{B}^+}
\newcommand{\btplus}{\widetilde{\mathbf{B}}^+}
\newcommand{\bt}{\widetilde{\mathbf{B}}}

\newcommand{\aplus}{\mathbf{A}^+}
\newcommand{\atplus}{\widetilde{\mathbf{A}}^+}
\newcommand{\at}{\widetilde{\mathbf{A}}}
\newcommand{\bhol}[1]{\mathbf{B}^+_{\mathrm{rig} #1}}

\newcommand{\e}{\mathbf{E}}
\newcommand{\eplus}{\mathbf{E}^+}
\newcommand{\etplus}{\widetilde{\mathbf{E}}^+}
\newcommand{\et}{\widetilde{\mathbf{E}}}

\newcommand{\dcris}{\mathbf{D}_{\mathrm{cris}}}
\newcommand{\ddR}{\mathbf{D}_{\mathrm{dR}}}
\newcommand{\dfont}{\mathbf{D}}

\newcommand{\dsen}{\mathbf{D}_{\mathrm{Sen}}}

\newcommand{\vale}{v_\mathbf{E}}

\newcommand{\ndr}{\mathbf{N}_{\mathrm{dR}}}

\title{Limites de repr\'esentations cristallines}

\author{Laurent Berger}
\address{MS 050 Brandeis University \\ PO Box 549110 \\
Waltham MA 02454-9110 \\ USA}
\email{laurent@brandeis.edu}
\urladdr{people.brandeis.edu/\~{}laurent}
\date{14 Mai 2002}

\begin{document}

\frontmatter

\begin{abstract}
On \'etudie quelques propri\'et\'es des $(\phi,\Gamma)$-modules 
associ\'es aux repr\'esentations absolument cristallines. Comme corollaire, 
on r\'epond (dans le cas ``non-ramifi\'e'') \`a deux questions de Fontaine. 
Premi\`erement, on montre qu'une $\Zp$-repr\'esentation, limite de
$\Zp$-repr\'esentations cristallines \`a poids de Hodge-Tate
born\'es est elle-m\^eme cristalline. Deuxi\`emement, on montre que
tout $\phi$-module filtr\'e admissible peut \^etre construit \`a partir d'un
$(\phi,\Gamma_F)$-module de $q$-hauteur finie 
(i.e. le foncteur $i^*: \mathbf{Is\boldsymbol{\Gamma}\boldsymbol{\Phi}M_S^+} \ra
\mathbf{MF_F^{ad,+}}$ est essentiellement surjectif). L'ingr\'edient principal est le calcul d'une
borne  effective pour l'annulateur
du conoyau de l'inclusion 
$\bplus \otimes_{\bplus_F} \dfont^+(V) \ra \bplus \otimes_{\Qp} V$.
\end{abstract}

\begin{altabstract}
We establish some properties of $(\phi,\Gamma)$-modules associated to absolutely crystalline 
representations. As a corollary, we can answer (in the ``unramified case'') two questions of
Fontaine. First, we show that a $\Zp$-representation, which is a limit of crystalline
$\Zp$-representations with bounded Hodge-Tate weights is itself crystalline. Second, we show
that every admissible filtered $\phi$-module can be constructed from a
$(\phi,\Gamma_F)$-module of finite $q$-height (that is, the functor 
$i^*: \mathbf{Is\boldsymbol{\Gamma}\boldsymbol{\Phi}M_S^+} \ra
\mathbf{MF_F^{ad,+}}$ is essentially surjective). The main ingredient is the computation of
an explicit bound for the annihilator of the cokernel of the inclusion 
$\bplus \otimes_{\bplus_F} \dfont^+(V) \ra \bplus \otimes_{\Qp} V$.
\end{altabstract}


\maketitle
\setcounter{tocdepth}{2}
\tableofcontents

\mainmatter

\section*{Introduction}
Soit $k$ un corps parfait de caract{\'e}ristique $p$,
$F$ le corps des fractions de l'anneau des vecteurs de Witt sur $k$
et $K$ une extension finie totalement ramifi{\'e}e de $F$;
la motivation de ce texte est de d\'emontrer une conjecture formul\'ee par Fontaine
\cite[Final Remark (c)]{Fo96}, dans le cas ``non-ramifi\'e'':
\begin{conj}\label{introconj}
Si $T$ est une $\Zp$-repr\'esentation de $G_F$, telle qu'il existe $b \geq a \in \NN$ et une suite
$\{T_i\}_i$ de $\Zp$-repr\'esentations cristallines de $G_F$ \`a poids de Hodge-Tate dans $[a;b]$
v\'erifiant $T/p^i \simeq T_i/p^i$, alors $T$ est cristalline.
\end{conj}

Si $b-a \leq p-1$, alors cette conjecture est un th\'eor\`eme, cons\'equence directe 
des constructions de Fontaine-Laffaille \cite{FL82}. L'analogue semi-stable de 
cette conjecture a \'et\'e d\'emontr\'e par Breuil 
dans \cite{Br}, dans le cas $b-a \leq p-2$.
On sait aussi faire les cas $b-a=0$ (pour une repr\'esentation de $G_K$, o\`u $K$ est
arbitrairement ramifi\'e), et $b-a=1$ (dans la plupart des cas, en utilisant 
les r\'esultats d\'emontr\'es par Breuil dans \cite{Br00}).

Notre d\'emonstration de la conjecture ci-dessus se fait en deux
\'etapes (toutes les notations employ\'ees ici sont d\'efinies dans le premier chapitre,
contentons-nous de mentionner que $\bplus_F=F \otimes_{\OO_F} \OO_F[[\pi]]$, et
que $q$ et les $\phi^{n-1}(q)$ 
sont des \'el\'ements de $\bplus_F$ d\'efinis par $\phi^{n-1}(q)=\Phi_{p^n}(1+\pi)$ o\`u
$\Phi_k$ est le $k$-i\`eme polyn\^ome cyclotomique.
En particulier, quotienter un $\bplus_F$-module par $\phi^{n-1}(q)$ revient \`a le ``localiser'' en
$\eps^{(n)}-1$). 

La th\'eorie des $(\phi,\Gamma)$-modules de Fontaine \cite{Fo91} et un th\'eor\`eme de Colmez
\cite{Co99} permettent d'associer \`a toute repr\'esentation cristalline $V$ de $G_F$ un
$\bplus_F$-module $\dfont^+(V)=(\bplus \otimes_{\Qp} V)^{H_F}$
libre de rang $d=\dim_{\Qp}(V)$, muni d'un Frobenius $\phi$ et d'une action de
$\Gamma_F=\on{Gal}(F(\mu_{p^{\infty}})/F)$ (on dit que $V$ est de hauteur finie). 
L'anneau $\bplus$ est un anneau de p\'eriodes construit
gr\^ace \`a la th\'eorie du corps de normes. 
La premi\`ere \'etape est de contr\^oler pr\'ecis\'ement
l'annulateur du conoyau de l'inclusion
$\bplus \otimes_{\bplus_F} \dfont^+(V) \ra \bplus \otimes_{\Qp} V$.

Ensuite Wach a montr\'e dans \cite{Wa96} qu'une repr\'esentation 
$V$ de hauteur finie est cristalline
si et seulement si $\dfont^+(V)$ contient un sous $\bplus_F$-module $N_V$ 
stable par $\Gamma_F$ et tel que $\Gamma_F$ agisse trivialement sur $N_V / \pi N_V$. 
On peut d'ailleurs imposer une l\'eg\`ere condition suppl\'ementaire sur $N_V$, et ce
$\bplus_F$-module est alors canoniquement associ\'e \`a $V$. La deuxi\`eme
\'etape est de construire \`a la main un tel $\bplus_F$-module pour une limite de repr\'esentations
cristallines.

Pour r\'ealiser la premi\`ere \'etape, nous sommes amen\'es \`a \'etudier de pr\`es les
$(\phi,\Gamma_F)$-modules associ\'es aux repr\'esentations absolument cristallines. Nous
montrons notamment:
\begin{theo}\label{introtheo}
Si $V$ est une repr\'esentation cristalline de $G_F$, 
dont les poids de Hodge-Tate sont $h_1 \leq \cdots \leq h_d$.
alors:
\begin{enumerate}
\item $\pi^{h_d-h_1} \bplus \otimes_{\Qp} V \subset \bplus \otimes_{\bplus_F} \dfont^+(V)$;
\item si $n\geq 1$ et $N_V$ est le $\bplus_F$-module construit par
Wach, alors l'action de $\Gamma_n=\on{Gal}(F_{\infty}/F_n)$ 
sur $V_n=N_V/\phi^{n-1}(q)$ 
est diagonale dans une base convenable de ce
$F_n$-espace vectoriel et $\gamma \in \Gamma_n$ agit par
$\on{Diag}(\chi(\gamma)^{h_1},\cdots,\chi(\gamma)^{h_d})$.
\end{enumerate}
\end{theo}
Le dernier point du th\'eor\`eme reste d'ailleurs vrai pour 
l'action de $\Gamma_n$ sur $\dfont^+(V)/\phi^{n-1}(q)$, car l'application naturelle
$N_V/\phi^{n-1}(q) \ra \dfont^+(V)/\phi^{n-1}(q)$ est un isomorphisme.

Pour r\'ealiser la deuxi\`eme \'etape, on montre la proposition suivante, qui nous permet de 
montrer la conjecture \ref{introconj} en construisant explicitement $N_T$ pour $T$.
\begin{prop}
Si $T_1$ et $T_2$ sont deux r\'eseaux de deux repr\'esentations cristallines \`a poids de
Hodge-Tate dans $[a;b]$, tels que $T_1 / p^n = T_2 / p^n$
avec $n \in \NN$ v\'erifiant $p^{n-1}(p-1) \geq b-a+1$, alors
$N_{T_1}$ et $N_{T_2}$ ont m\^eme image dans $(\aplus/p^n \otimes_{\Zp/p^n} T_i/p^n)^{H_F}$.
\end{prop}
Cette proposition montre que si les $T_i$ varient continument, alors il en va de m\^eme
pour les $N_{T_i}$, sur un disque dont le rayon diminue quand on augmente la longueur de la
filtration.

Le th\'eor\`eme \ref{introtheo} (plus un petit calcul) nous permet de
d\'emontrer une autre conjecture de Fontaine. 
Soient $\mathbf{Is\boldsymbol{\Gamma}\boldsymbol{\Phi}M_S^+}$ la cat\'egorie 
des $\bplus_F$-modules libres $N$ de rang fini munis d'un Frobenius $\phi$
tel que $N/\phi^*N$ est annul\'e
par une puissance de $q$, et $\mathbf{MF_F^{ad,+}}$ la cat\'egorie des
$\phi$-modules filtr\'es admissibles. On peut munir tout objet $N$ de
$\mathbf{Is\boldsymbol{\Gamma}\boldsymbol{\Phi}M_S^+}$
d'une filtration qui d\'epend de $\phi$, et $i^*N=N/\pi N$ est alors un $\phi$-module filtr\'e.
Fontaine a conjectur\'e dans \cite[2.3.6]{Fo91} que le foncteur $i^*: 
\mathbf{Is\boldsymbol{\Gamma}\boldsymbol{\Phi}M_S^+} \ra \mathbf{MF_F^{ad,+}}$
est essentiellement surjectif. De fait:
\begin{theo}
Si $V$ est une repr\'esentation cristalline positive de $G_F$, 
et si $N_V$ est le $\bplus_F$-module construit par Wach,  
muni de la filtration $\on{Fil}^i N_V = \{ x \in N_V, \phi(x) \in
q^i N_V \}$, alors $\dcris(V) \subset \bhol{,F} \otimes_{\bplus_F} N_V$ et
l'application naturelle 
$\dcris(V) \ra N_V / \pi N_V$ est 
un isomorphisme de $\phi$-modules filtr\'es.
\end{theo}
On en d\'eduit imm\'ediatement la conjecture mentionn\'ee ci-dessus. 

On termine avec une nouvelle d\'emonstration de certains des r\'esultats de
\cite{FL82,Wa97} (construction de r\'eseaux galoisiens associ\'es \`a des r\'eseaux fortement
divisibles, quand la longueur de la filtration est $\leq p-2$).

\begin{rema}
Dans une version ant\'erieure de cet article, je d\'emontrais  que $N_T/\pi N_T$ s'identi\-fiait
\`a un r\'eseau fortement divisible de $\dcris(V)$. Ceci est incorrect. Par exemple,
comme me l'ont fait remarquer C. Breuil, F. Diamond et P. Colmez,
certaines repr\'esentations attach\'ees \`a des formes modulaires sont irr\'eductibles, mais
pas leur r\'eduction modulo $p$, et dans certains cas on obtient plus de r\'eseaux
galoisiens que de modules fortement divisibles, ce qui aboutit \`a une contradiction.
\end{rema}

\noindent\textbf{\itshape Remerciements:} Je souhaite remercier Pierre Colmez, qui est \`a
l'origine de plusieurs des id\'ees de cet article; il a aussi relu 
attentivement de nombreuses versions
de ce texte. Lors de la r\'edaction de cet article, j'ai b\'en\'efici\'e de conversations
enrichissantes avec lui et avec Christophe Breuil, Fred Diamond, Jean-Marc Fontaine et Barry Mazur. 
Enfin, l'article ``Repr\'esentations
$p$-adiques potentiellement cristallines'' de Nathalie Wach a \'et\'e une source constante
d'inspiration.

\section{Rappels et notations}
Ce chapitre est consacr\'e \`a quelques rappels sur les p\'eriodes $p$-adiques; on pourra se
reporter \`a \cite{Moi} pour beacoup des constructions d\'ecrites ci-dessous.
Soit $k$ un corps parfait de caract{\'e}ristique $p$,
$F$ le corps des fractions de l'anneau des vecteurs de Witt sur $k$
et $K$ une extension finie totalement ramifi{\'e}e de $F$. Soit
$\overline{F}$ une
cl{\^o}ture alg{\'e}brique de $F$ et $\Cp=\hat{\overline{F}}$ sa compl{\'e}tion
$p$-adique. 
On pose $G_K=\on{Gal}(\overline{K}/K)$, c'est aussi le groupe
des automorphismes continus $K$-lin{\'e}aires de $\Cp$.
Le corps $\Cp$ est un corps complet alg{\'e}briquement clos
de corps r{\'e}siduel $\overline{k}$. 
On pose aussi $K_n=K(\mu _{p^n})$ et $K_{\infty}$ est d{\'e}fini
comme  {\'e}tant la r{\'e}union des $K_n$.  
Soit $H_K$ le noyau du caract{\`e}re cyclotomique 
$\chi: G_K \ra \Zp^*$ et $\Gamma_K=G_K/H_K$
le groupe de Galois de $K_{\infty}/K$
qui s'identifie via le caract{\`e}re cyclotomique {\`a} un sous groupe ouvert
de $\Zp^*$. Si $K=F$, on a d'ailleurs $\Gamma_F \simeq \Zp^*$, et c'est un groupe pro-cyclique
sauf si $p=2$; si $p \neq 2$ ou $p=2$ et $n \geq 2$, alors
$\Gamma_n=\on{Gal}(K_{\infty}/K_n)$ est procyclique.

Dans la suite, $T$ d\'esignera
une $\Zp$-repr\'esentation de $G_K$, c'est \`a dire que $T$ est un $\Zp$-module libre de
rang $d$, muni d'une action lin\'eaire et continue de $G_K$, et $V=\Qp \otimes_{\Zp} T$.
La principale strat{\'e}gie 
(due {\`a} Fontaine, voir par exemple \cite{Bu88sst})
pour {\'e}tudier les
repr{\'e}sentations $p$-adiques d'un groupe $G$ est de cons\-truire des 
$\Qp$-alg{\`e}bres topologiques $B$ munies d'une action du groupe
$G$ et de structures suppl{\'e}\-mentaires de telle mani{\`e}re que si $V$
est une repr{\'e}sentation $p$-adique, 
alors $D_B(V)=(B \otimes_{\Qp} V)^G$  est un
$B^G$-module qui h{\'e}rite de ces structures, et que le foncteur
qui {\`a} $V$ associe $D_B(V)$ fournisse des invariants int{\'e}ressants
de $V$. On dit qu'une repr{\'e}sen\-tation $p$-adique $V$ de $G$ est
\emph{$B$-admissible} si on a $B\otimes_{\Qp} V \simeq B^d$ en tant que
$G$-modules. 

\subsection{Le corps $\et$ et ses sous-anneaux}
Soient \[ \et=\projlim_{x\mapsto x^p} \Cp 
=\{ (x^{(0)},x^{(1)},\cdots) \mid (x^{(i+1)})^p = x^{(i)} \} \]
et $\etplus$ l'ensemble des $x \in \et$ tels que $x^{(0)} \in \OO_{\Cp}$.
Si $x=(x^{(i)})$ et $y=(y^{(i)})$ sont deux {\'e}l{\'e}ments de $\et$,
alors on d{\'e}finit leur somme $x+y$ et leur produit $xy$ par:
\[ (x+y)^{(i)}= \lim_{j \ra \infty} (x^{(i+j)}+y^{(i+j)})^{p^j} \text{ et }
(xy)^{(i)}=x^{(i)}y^{(i)} \] ce qui fait de
$\et$ un corps de caract{\'e}ristique $p$ dont on peut montrer qu'il est
alg{\'e}briquement clos.
Si $x=(x^{(n)})_{n \geq 0} \in \et$ soit $\vale(x)=v_p(x^{(0)})$. C'est
une valuation sur $\et$ pour laquelle celui-ci est complet;
l'anneau des entiers de $\et$ est $\etplus$.
Soit $\atplus$ l'anneau $W(\etplus)$ des vecteurs de Witt {\`a}
coefficients dans $\etplus$ et $\btplus =\atplus[1/p] 
=\{ \sum_{k\gg -\infty} p^k [x_k],\ x_k \in \etplus
\}$ 
o{\`u} $[x] \in \atplus$ est le
rel{\`e}vement de Teichm{\"u}ller de $x \in \etplus$.
Cet anneau est muni d'un morphisme d'anneaux
$\theta: \btplus \ra \Cp$ d{\'e}fini par la formule
$\theta(\sum_{k\gg -\infty} p^k[x_k] )
=\sum_{k\gg -\infty} p^k x_k^{(0)}$.
Soient $\eps=(\eps^{(i)})\in\etplus$ avec $\eps^{(0)}=1$ 
et $\eps^{(1)}\neq 1$, $\pi=[\eps]-1$, $\pi_1=
[\eps^{p^{-1}}]-1$, $\omega=\pi/\pi_1$
et $q=\phi(\omega)=\phi(\pi)/\pi$. 
Alors $\ker(\theta)$
est l'id{\'e}al principal engendr{\'e} par $\omega$.

Remarquons que $\eps$ est un {\'e}l{\'e}ment de $\etplus$ tel que
$\vale(\eps-1)=p/(p-1)$. On pose $\e_F=k((\eps-1))$ et on d{\'e}finit
$\e$ comme {\'e}tant 
la cl{\^o}ture s{\'e}parable de $\e_F$ dans $\et$ ainsi que $\eplus=\e
\cap \etplus$ l'anneau des entiers de $\e$. 
Remarquons que, par d{\'e}finition, $\e$ est s{\'e}parablement clos, et que 
l'on retrouve $\et$ {\`a} partir de $\e$ en prenant le compl{\'e}t{\'e} de
sa cl{\^o}ture radicielle.

\subsection{L'anneau $\bdR$ et ses sous-anneaux}
L'anneau $\bdR^+$ est d{\'e}fini comme {\'e}tant le compl{\'e}t{\'e} pour la 
topologie $\ker(\theta)$-adique de $\btplus$ (on remarquera que
$\atplus$ est complet pour cette topologie):
\[ \bdR^+=\projlim_{n\geq 0} \btplus/(\ker(\theta)^n) \]
c'est un anneau de valuation discr{\`e}te, d'id{\'e}al maximal $\ker(\theta)$;
la s{\'e}rie qui d{\'e}finit $\log([\eps])$ converge
dans $\bdR^+$ vers un {\'e}l{\'e}ment $t$, qui est un g{\'e}n{\'e}rateur de l'id{\'e}al
maximal, ce qui fait que $\bdR=\bdR^+[1/t]$ est un corps, muni d'une
action de $G_F$ et d'une filtration d{\'e}finie par 
$\on{Fil}^i(\bdR)=t^i \bdR^+$ pour $i \in \ZZ$.

On dit qu'une repr{\'e}sentation $V$ de $G_K$
est de de Rham si 
elle est $\bdR$-admissible ce qui {\'e}quivaut {\`a} ce que
le $K$-espace
vectoriel $\ddR(V)=(\bdR\otimes_{\Qp} V)^{G_K}$
est de dimension $d=\dim_{\Qp}(V)$.

L'anneau $\bmax^+$ est d{\'e}fini comme {\'e}tant
\[ \bmax^+= \{ \sum_{n \geq 0} a_n \frac{\omega^n}{p^n}
\text{ o{\`u} $a_n\in \btplus$ est une suite qui tend vers $0$} \} \]
et $\bmax=\bmax^+[1/t]$. 
On peut d'ailleurs remplacer $\omega$ par
n'importe quel g{\'e}n{\'e}rateur de $\ker(\theta)$.
Cet anneau se plonge canoniquement dans $\bdR$ (les s{\'e}ries
d{\'e}finissant ses {\'e}l{\'e}ments convergent dans $\bdR$) 
et en particulier il est muni de l'action
de Galois et de la filtration induites par celles de $\bdR$, ainsi que
d'un Frobenius $\phi$, qui {\'e}tend l'application $\phi:\atplus \ra
\atplus$ d{\'e}duite de $x \mapsto x^p$ dans $\etplus$.
On remarquera que $\phi$ ne se prolonge pas par continuit{\'e} {\`a} $\bdR$.
On pose aussi $\btrigplus{}=\cap_{n=0}^{+\infty}
\phi^n(\bmax^+)$.

On dit qu'une repr{\'e}sentation $V$ de $G_K$
est cristalline si 
elle est $\bmax$-admissible ou, ce qui revient au m{\^e}me,
$\btrigplus{}[1/t]$-admissible
(les p{\'e}riodes des repr{\'e}sentations cristallines vivent dans des
sous $F$-espaces vectoriels de dimension finie,
stables par $\phi$, 
de $\bmax$, et donc en fait
dans $\cap_{n=0}^{+\infty}
\phi^n(\bmax^+)[1/t]$); 
ceci {\'e}quivaut {\`a} ce que
le $F$-espace
vectoriel \[ \dcris(V)=(\bmax \otimes_{\Qp} V)^{G_K}=(\btrigplus{}[1/t]
\otimes_{\Qp} V)^{G_K} \] 
est de dimension 
$d=\dim_{\Qp}(V)$. Alors $\dcris(V)$ est muni d'un Frobenius et d'une 
filtration induits par ceux de $\bmax$, 
et $(\bdR\otimes_{\Qp} V)^{G_K}=\ddR(V)=
K \otimes_F \dcris(V)$ ce qui fait qu'une repr{\'e}sentation cristalline
est aussi de de Rham.

Si $V$ est une repr\'esentation $p$-adique, on dit que $V$ est de Hodge-Tate, \`a poids de
Hodge-Tate $h_1, \cdots, h_d$, si l'on a une d\'ecomposition $\Cp\otimes_{\Qp} V \simeq
\oplus_{j=1}^d \Cp(h_j)$. Dans ce cas, on voit que
$(\Cp\otimes_{\Qp} V)^{H_K} \simeq \oplus_{j=1}^d \hat{K}_{\infty} (h_j)$ et on peut montrer que la
r\'eunion $\dsen(V)=(\Cp\otimes_{\Qp} V)^{H_K}_{\mathrm{fini}}$
des sous $K_{\infty}$-espaces vectoriels de
dimension finie stables par $\Gamma_K$ de $(\Cp\otimes_{\Qp} V)^{H_K}$ est \'egale \`a
$\oplus_{j=1}^d K_{\infty} (h_j)$. Le $K_{\infty}$-espace vectoriel 
$\dsen(V)$ est muni d'une action r\'esiduelle de $\Gamma_K$, et si $\gamma \in \Gamma_K$ est un
\'el\'ement  suffisamment proche de $1$, alors la s\'erie d'op\'erateurs
$\log(\gamma)/\log_p(\chi(\gamma))$ converge vers un op\'erateur $K_{\infty}$-lin\'eaire 
$\nabla_V : \dsen(V) \ra \dsen(V)$ 
qui ne d\'epend pas du choix de $\gamma$, et 
qui est diagonalisable \`a valeurs propres $h_1, \cdots,
h_d$.
On dira qu'une repr\'esentation $p$-adique 
$V$ est positive si ses poids de Hodge-Tate sont $\leq 0$
(la d\'efinition du signe des poids de Hodge-Tate est malheureuse).

\subsection{L'anneau $\bt$ et ses sous-anneaux}
Soit $\at$ l'anneau des vecteurs de Witt {\`a} coefficients dans
$\et$ et $\bt=\at[1/p]$.
Soit $\mathbf{A}_F$
le compl{\'e}t{\'e} de $\OO_F[\pi,\pi^{-1}]$ dans
$\at$ pour la topologie de celui-ci, c'est aussi le compl{\'e}t{\'e}
$p$-adique de $\OO_F[[\pi]][\pi^{-1}]$. 
C'est un anneau de valuation discr{\`e}te complet
dont le corps r{\'e}siduel est $\mathbf{E}_F$.
Soit $\mathbf{B}$ le compl{\'e}t{\'e}
pour la topologie $p$-adique
de l'extension maximale non ramifi{\'e}e de 
$\mathbf{B}_F=\mathbf{A}_F[1/p]$
dans $\bt$. On d{\'e}finit alors 
$\mathbf{A}=\mathbf{B}
\cap \at$
et $\aplus=\mathbf{A} \cap \atplus$. Ces anneaux sont 
munis d'une action de Galois et d'un Frobenius d{\'e}duits de ceux de $\et$. 
On pose $\mathbf{A}_K=\mathbf{A}^{H_K}$ et
$\mathbf{B}_K=\mathbf{A}_K[1/p]$.
Quand $K=F$ les deux d{\'e}finitions co{\"\i}ncident.
On pose de m\^eme $\bplus=\aplus[1/p]$ et $\aplus_F=(\aplus)^{H_F}$ ainsi que
$\bplus_F=(\bplus)^{H_F}$.

Si $V$ est une repr{\'e}sentation $p$-adique de $G_K$ soit $\dfont(V)
=(\mathbf{B} \otimes_{\Qp} V)^{H_K}$. 
On sait \cite{Fo91} que $\dfont(V)$ est un
$\mathbf{B}_K$-espace vectoriel de dimension $d=\dim(V)$ muni d'un Frobenius
et d'une action r{\'e}siduelle de $\Gamma_K$ 
qui commutent (c'est un $(\phi,\Gamma_K)$
module) et que l'on peut r{\'e}cup{\'e}rer $V$ gr{\^a}ce {\`a} la formule
$V=(\mathbf{B} \otimes_{\mathbf{B}_K} \dfont(V))^{\phi=1}$. 

On se place maintenant dans le cas $K=F$.
On dit qu'une repr{\'e}sentation $p$-adique $V$ de $G_F$ 
est de hauteur finie si
$\dfont(V)$ poss{\`e}de une base sur $\mathbf{B}_F$ form{\'e}e
d'{\'e}l{\'e}ments de $\dfont^+(V) = (\bplus \otimes_{\Qp} V)^{H_F}$. Un
r{\'e}sultat de Fontaine 
\cite{Fo91} (voir aussi \cite[III.2]{Co99}) montre que
$V$ est de hauteur finie si et
seulement si $\dfont(V)$ poss{\`e}de un sous-$\bplus_F$-module libre de
type fini stable par $\phi$
de rang {\'e}gal {\`a} $d=\dim_{\Qp}(V)$.
Rappelons le r\'esultat principal (le th\'eor\`eme de Colmez \cite{Co99}, voir aussi
\cite{Moi}) sur les repr\'esentation cristallines de $G_F$:
\begin{theo}\label{hf}
Si $V$ est une repr{\'e}sentation cristalline de $G_F$, 
alors $V$ est de hauteur finie.
\end{theo}

\subsection{Compl\'ements sur $\bplus_F$ et $\bhol{,F}$}
Soit $\bhol{,F}$ l'ensemble des s\'eries de la forme 
$A(\pi)=\sum_{k=0}^{\infty} a_k \pi^k$,
telles que $a_k \in F$ et telles que la s{\'e}rie  $A(T)$
converge sur le disque unit\'e ouvert
$\{ z \in \Cp, |z| < 1 \}$. Rappelons que $\bplus_F$ est un anneau
principal, et que $\bhol{,F}$ se comporte comme un anneau principal, en ce sens
qu'il est de Bezout (tout id\'eal de type fini est principal) et 
\cite{Moi} qu'il admet la th\'eorie
des diviseurs \'el\'ementaires.

Si $R$ est un anneau qui a la th\'eorie
des diviseurs \'el\'ementaires, 
et si $M$ est un sous-module de type fini de
$N \simeq R^d$, alors il existe des \'el\'ements $r_1, \cdots, r_d$ de $R$ et une base
$n_1, \cdots, n_d$ de $N$ tels que $r_1 | \cdots | r_d$  et $r_1 n_1, \cdots r_d n_d$
est une base de $M$, et les id\'eaux $(r_1), \cdots, (r_d)$ sont d\'etermin\'es
de mani\`ere unique 
par ces conditions: c'est le th\'eor\`eme des diviseurs \'el\'ementaires.
On \'ecrira pour simplifier: $[N:M]=[r_1;\cdots;r_d]$.

Les deux anneaux $\bplus_F$ et $\bhol{,F}$ sont munis d'actions de $\phi$ et de $\Gamma_F$ qui
commutent, et que l'on d\'eduit par semi-lin\'earit\'e des formules $\phi(T)=(1+T)^p-1$ et
$\gamma(T)=(1+T)^{\chi(\gamma)}-1$.

Si $n \geq 1$, alors l'id\'eal
engendr\'e par $\phi^{n-1}(q)$ est stable par
$\Gamma_F$ et on a une identification de $F[\Gamma_F]$-modules
$\bhol{,F}/\phi^{n-1}(q)=\bplus_F/\phi^{n-1}(q)=F_n$.

\begin{lemm}\label{idealphi}
Si $I$ est un id\'eal principal de $\bhol{,F}$, qui est stable par $\Gamma_F$, alors
$I$ est engendr\'e par un \'el\'ement de la forme
$\pi^{j_0}\prod_{n=1}^{+\infty}(\phi^{n-1}(q)/p)^{j_n}$. De plus:
\begin{enumerate}
\item $\phi(I) \subset I$ si et seulement si la suite $\{j_n\}_n$ est
d\'ecroissante;
\item $I \subset \phi(I)$ si et seulement si la suite $\{j_n\}_n$ est
croissante;
\item les r\'esultats pr\'ec\'edents restent vrais pour les id\'eaux de $\bplus_F$, et dans
ce cas $j_n=0$ si $n \gg 0$ (par exemple, $I \subset \phi(I)$ implique que $I=\bplus_F$).
\end{enumerate}
\end{lemm}

\begin{proof}
Voir \cite[III.8]{Co99} pour une d\'emonstration du fait que
si $I$ est un id\'eal de $\bplus_F$, qui est stable par $\Gamma_F$, alors
$I$ est engendr\'e par un \'el\'ement de la forme
$\pi^{j_0}\prod_{n=1}^{n_0}(\phi^{n-1}(q)/p)^{j_n}$. Cette d\'emonstration s'adapte imm\'ediatement au
cas de $\bhol{,F}$. Enfin, les trois points qui suivent ne posent pas de difficult\'e.
\end{proof}

\begin{exem}
L'id\'eal de $\bhol{,F}$ engendr\'e par $t=\log(1+\pi)$ satisfait les deux
points ci-dessus. Cela correspond \`a la d\'ecomposition $t=\pi
\prod_{n=1}^{+\infty}(\phi^{n-1}(q)/p)$.
\end{exem}

\section{Une bonne description de $\dcris(V)$}
Dans toute la suite, $V$ d\'esignera une repr\'esentation cristalline de $G_F$, 
de dimension $d$, dont les
poids de Hodge-Tate sont n\'egatifs. On notera $r_1 \leq \cdots \leq r_d$ les oppos\'es des
poids de Hodge-Tate de $V$.

\subsection{La construction de Wach}\label{constnw}
L'objet de ce paragraphe est de pr\'eciser les r\'esultats de Wach sur les repr\'esentations
cristallines de hauteur finie. Rappelons que Wach a d\'emontr\'e dans \cite[p. 380]{Wa96} le
th\'eor\`eme suivant:
\begin{theo}
Si $V$ est une repr\'esentation de hauteur finie de $G_F$, alors $V$ est cristalline positive si et
seulement s'il existe un sous $\bplus_F$-module $N$ libre de rang $d$ de $\dfont^+(V)$ stable par
$\Gamma_F$ et tel que l'action de $\Gamma_F$ soit triviale sur $N/\pi N$.
\end{theo}
Si $V$ est une repr\'esentation cristalline positive de $G_F$, alors le th\'eor\`eme de Colmez
montre que $V$ est de hauteur finie. Il existe donc $N \subset \dfont^+(V)$ tel que
$N$ est un $\bplus_F$-module libre de rang $d$ stable par
$\Gamma_F$ et tel que l'action de $\Gamma_F$ soit triviale sur $N/\pi N$.
\begin{lemm}
Sous les hypoth\`eses pr\'ec\'edentes, on peut s'arranger pour qu'il existe un entier $s$ tel que
$\pi^s \dfont^+(V) \subset N \subset \dfont^+(V)$. Un tel $N$ est alors unique
(en particulier il est stable par $\phi$), et on le notera $N_V$.
\end{lemm}
\begin{proof}
L'id\'eal $I_N$ de
$\bplus_F$ constitu\'e de l'ensemble des $\lambda \in \bplus_F$ tels que $\lambda \dfont^+(V)
\subset N$  est stable par $\Gamma_F$, et il est donc engendr\'e par un \'el\'ement 
de la forme $\pi^{\alpha_0} q^{\alpha_1} \cdots (\phi^{s-1}(q))^{\alpha_s}$. 
Si $N'=\dfont^+(V) \cap N[\phi^{n-1}(q^{-1})]_{n \geq 1}$, alors $N'$ 
est un sous-module libre de rang $d$ de
$\dfont^+(V)$, stable par $\Gamma_F$,
et l'application $N /\pi N \ra N' /\pi N'$ est un isomorphisme (car elle est injective)
ce qui fait que $\Gamma_F$ agit trivialement sur $N' / \pi N'$. Enfin, on voit que si 
$x \in \dfont^+(V)$ est tel que $\pi^{\beta_0} q^{\beta_1} \cdots (\phi^{s-1}(q))^{\beta_s} x
\in N'$, alors $\pi^{\beta_0} x \in N'$, ce qui fait que $\pi^s \dfont^+(V) \subset N'$ avec
$s=\alpha_0$.

Montrons qu'un tel $N$ est unique. S'il existait $N_1$ et $N_2$ satisfaisant les conditions du lemme,
alors pour tout $x \in N_1 \setminus \pi N_1$, 
il existerait un unique $s \in \NN$ tel que $\pi^s x \in N_2 \setminus \pi
N_2$. Comme l'action de $\Gamma_F$ est triviale sur $N_i / \pi N_i$, on voit que $s=0$. En appliquant
cela \`a une base de $N_1$, on voit que $N_1 \subset N_2$ et donc par sym\'etrie que $N_1=N_2$.

Enfin, on voit que $N+\phi^*N$ satisfait les m\^emes hypoth\`eses que $N$, et donc $N+\phi^*N=N$ ce
qui fait que $N$ est stable par $\phi$.
\end{proof}

Si $T$ est un r\'eseau de $V$, on pose $\dfont^+(T)=(\aplus \otimes_{\Zp} T)^{H_F}$, c'est un
$\aplus_F$-r\'eseau de $\dfont^+(V)$, de m\^eme que
$\dfont(T)=(\mathbf{A} \otimes_{\Zp} T)^{H_F}$ est un
$\mathbf{A}_F$-r\'eseau de $\dfont(V)$
\begin{lemm}
Si $T$ est un r\'eseau de $V$, et $N_T = N_V \cap \dfont^+(T)$, alors $N_T$ est un sous
$\aplus_F$-module libre de rang $d$ de $\dfont^+(T)$, stable par $\phi$ et $\Gamma_F$  tel que
l'action de $\Gamma_F$ sur $N_T / \pi N_T$ soit triviale, et tel 
qu'il existe un entier $s \geq 0$ tel que $\pi^s \dfont^+(T) \subset N_T$. Ceci d\'efinit $N_T$ de
mani\`ere unique.
\end{lemm}
\begin{proof}
La seule chose qui n'est pas \'evidente est que $N_T$ est libre (puisque $\aplus_F$ n'est pas
principal). Mais on voit que $N_T = N_V \cap \dfont(T)$, et $N_T$ est alors libre par
\cite[B.1.2.4]{Fo91}.
\end{proof}
On a donc d\'emontr\'e la proposition suivante:
\begin{prop}
Si $T$ est un r\'eseau d'une repr\'esentation cristalline positive $V$, alors il existe un unique
sous $\aplus_F$-module $N_T$
libre de rang $d$ de $\dfont^+(T)$ tel que $N_T$ est stable par $\Gamma_F$ et
l'action de $\Gamma_F$ est triviale sur $N_T / \pi N_T$, et tel qu'il existe un entier $s \geq 0$ tel
que $\pi^s \dfont^+(T) \subset N_T$. Dans ce cas $N_T$ est stable par $\phi$, et le sous
$\bplus_F$-module $N_V =  \bplus_F \otimes_{\aplus_F} N_T$ de $\dfont^+(V)$ v\'erifie
les conditions correspondantes.
\end{prop}

\subsection{L'inclusion $\dcris(V) \subset \dhol(V)$}
Soit $\bnrig{}{,F}$ l'ensemble des s\'eries 
de la forme
$A(\pi)=\sum_{k=-\infty}^{+\infty} a_k \pi^k$,
telles que $a_k \in F$ et qu'il existe $r<1$ tel que la s{\'e}rie  $A(T)$
converge sur la couronne 
$\{ z \in \Cp, r < |z| < 1 \}$. 

On a montr\'e dans \cite{Moi} que si $V$ est une repr\'esentation cristalline 
positive de $G_F$ 
(et donc dans ce cas de hauteur finie), alors
$\dcris(V) \subset \dnrig{}(V)=\bnrig{}{,F}
\otimes_{\bplus_F} \dfont^+(V)$. 
On a en fait dans notre cas particulier un meilleur contr\^ole des
p\'eriodes de $V$; soit
$\dhol(V)=\bhol{,F}\otimes_{\bplus_F}\dfont^+(V)$, on a alors:
\begin{prop}\label{crishol}
Si $V$ est une repr\'esentation cristalline positive de $G_F$, 
alors $\dcris(V) \subset \bhol{,F} \otimes_{\bplus_F} N_V$. En particulier 
$\dcris(V) \subset \dhol(V)$, et $\dcris(V)=\dhol(V)^{\Gamma_F}=(\bhol{,F}
\otimes_{\bplus_F} N_V)^{\Gamma_F}$.
\end{prop}

\begin{proof}
Soit $\{e_i\}$ une base de
$N_V$, et $P$ et $G=G(\gamma)$ les matrices de $\phi$ et $\gamma \in \Gamma_F \setminus \{ 1\}$ 
dans cette base. Soit
$\{x_i\}$ une base de $\dcris(V)$, on sait \cite{Moi} que $x_i \in \dnrig{}(V)$, et donc si
$M$ est telle que $(x_i)=M(e_i)$, alors $M \in \on{M}(d,\bnrig{}{,F})$. Si $D \in \on{M}(d,F)$ est la
matrice de $\phi$ dans la base $\{x_i\}$, alors on a $\gamma(M)G=M$ et $\phi(M)P=DM$.

Montrons qu'il existe des entiers $\alpha_j$ et $u \in (\bplus_F)^*$ tels que $\delta=\det(P)=
q^{\alpha_1}\cdots (\phi^{s-1}(q))^{\alpha_s} u$.
Si $g$ est le d\'eterminant de $G$, alors $g=1\mod{\pi}$ et comme les actions de
$\phi$ et $\Gamma_F$ commutent, on a $\gamma(\delta)g=\phi(g)\delta$, ce qui fait  que $\delta$
et $\gamma(\delta)$ engendrent le m\^eme id\'eal de $\bplus_F$. Par le lemme 
\ref{idealphi}, c'est
donc que $\delta=\pi^{\alpha_0}q^{\alpha_1}\cdots (\phi^{s-1}(q))^{\alpha_s} u$ o\`u $u \in
(\bplus_F)^*$. De plus, si l'on regarde l'\'equation $\gamma(\delta)/\delta=\phi(g)/g$ modulo
$\pi$, on voit que $\alpha_0=0$. 

Posons $\nu=\pi^{\alpha_1} \cdots
\phi^{s-1}(\pi)^{\alpha_s}$, on a $\delta=\phi(\nu)\nu^{-1}u$.
Soient $N=M\nu$ et $Q=\delta P^{-1}$,
c'est la transpos\'ee de la matrice des cofacteurs de $P$, et donc $Q\in \on{M}(d,\bplus_F)$.
On a \[ \phi(N) = \phi(M) \phi(\nu) 
= DM \phi(\nu) P^{-1} 
= DM (\nu \delta u^{-1})  \delta^{-1}Q = DNu^{-1}Q \]
et l'argument de r\'egularisation de \cite{Moi} montre que si $\phi(N)=DNu^{-1}Q$, alors $N
\in \on{M}(d,\bhol{,F})$. Il existe donc $\lambda \neq 0$, $\lambda \in \bplus_F$, tel que  
$\lambda \dcris(V) \subset \bhol{,F} \otimes_{\bplus_F} N_V$.
 
Si $I$ est l'id\'eal de $\bplus_F$ form\'e
de l'ensemble des $\lambda \in \bplus_F$ tels que  
$\lambda \dcris(V) \subset \bhol{,F} \otimes_{\bplus_F} N_V$, alors $I$ est un id\'eal de
$\bplus_F$ qui est stable par $\Gamma_F$ et $\phi$, et 
par le lemme \ref{idealphi} il est engendr\'e par un \'el\'ement de la
forme $\lambda=\pi^{\beta_0} \cdots (\phi^{s-1}(q))^{\beta_s}$, avec $\beta_0 \geq \cdots \geq
\beta_s \geq 0$. 

Reste \`a montrer que $\lambda=1$, ce qui revient \`a montrer que
$\beta_0=0$ puisque $\beta_j \leq \beta_0$.
Rappelons que le groupe $\Gamma_F$ agit trivialement sur
$N_V /\pi N_V$, et donc sur $\bhol{,F} \otimes_{\bplus_F} N_V /\pi \bhol{,F}
\otimes_{\bplus_F} N_V$.  
Il existe $y\in\dcris(V)$ tel que $y\lambda \in \bhol{,F} \otimes_{\bplus_F} N_V  \setminus
\pi \bhol{,F} \otimes_{\bplus_F} N_V$. L'image de $y\lambda$ dans $\bhol{,F}
\otimes_{\bplus_F} N_V /\pi \bhol{,F} \otimes_{\bplus_F} N_V$  est non-nulle et $\Gamma_F$
agit dessus par $\chi^{\beta_0}$, ce qui fait que $\beta_0 = 0$.
\end{proof}

\begin{rema}
On voit d'ailleurs que 
$\bhol{,F}\otimes_F \dcris(V) \simeq \cap_{n=0}^{+\infty} \phi^{*n}
(\bhol{,F} \otimes_{\bplus_F} N_V)$, ce qui fournit une description 
assez simple de $\bhol{,F}\otimes_F \dcris(V)=\ndr(V)$, o\`u $\ndr(V)$ est le
module construit dans \cite{Moi}.
\end{rema}

Soient $(\lambda_i)$ les id\'eaux de $\bhol{,F}$ d\'efinis
par $[\bhol{,F} \otimes_{\bplus_F} N_V :\bhol{,F} \otimes_F
\dcris(V)]=[\lambda_1;\cdots;\lambda_d]$. 

\begin{prop}\label{invrigcris}
Il existe des entiers $\beta_{n,i}$ tels que
les id\'eaux $(\lambda_i)$ sont engendr\'es par 
$\lambda_i = \pi^{\beta_{0,i}} \prod_{n=1}^{\infty} (\phi^{n-1}(q) / p)^{\beta_{n,i}}$. 
\end{prop}

\begin{proof}
On a montr\'e au paragraphe pr\'ec\'edent que $\bhol{,F} \otimes_F
\dcris(V) \subset 
\bhol{,F} \otimes_{\bplus_F} N_V$,
et les deux modules en question sont munis d'une action de $\Gamma_F$, telle que $\gamma 
\in \Gamma_F$ agit par un isomorphisme.
Ceci montre (par unicit\'e des diviseurs \'el\'ementaires) que les id\'eaux
$(\lambda_i)$ et $\gamma(\lambda_i)$ sont \'egaux, et
le lemme \ref{idealphi} montre que l'on peut prendre
$\lambda_i=\pi^{\beta_{0,i}} \prod_{n=1}^{\infty} (\phi^{n-1}(q) / p)^{\beta_{n,i}}$, o\`u
$\{\beta_{n,i}\}_n$ est une suite d'entiers.
\end{proof}

\section{L'inclusion $\dfont^+(V) \subset \bplus \otimes_{\Qp} V$}
Dans toute la suite, on \'ecrira $\on{Diag}(a_1,\cdots,a_d)$ pour 
d\'esigner la matrice $d \times d$:
\[ \begin{pmatrix}
a_1 & & \\
 & \ddots & \\
 & & a_d
\end{pmatrix} \]
\subsection{L'action de $\Gamma_F$ sur $N_V/\phi^{n-1}(q)$}\label{gammaqn}
Si $n \geq 1$, on pose $\Gamma_n=\on{Gal}(F_{\infty}/F_n)$.
\begin{prop}
L'action de $\Gamma_n$ sur $V_n=N_V/\phi^{n-1}(q)$ 
est diagonale dans une base convenable de ce
$F_n$-espace vectoriel et $\gamma \in \Gamma_n$ agit dans cette base par
$\on{Diag}(\chi(\gamma)^{-\beta_{n,1}},\cdots,\chi(\gamma)^{-\beta_{n,d}})$.
\end{prop}

\begin{proof}
Comme $\Gamma_n$ est ab\'elien, il suffit de montrer la proposition pour un \'el\'ement
$\gamma_n$ de $\Gamma_n$, car une famille commutative d'endomorphismes diagonalisables est
codiagonalisable (si $p \neq 2$ ou $n \geq 2$, il suffit de prendre pour $\gamma_n$
un g\'en\'erateur topologique de $\Gamma_n$). 

Si $\Lambda=\on{Diag}(\lambda_1,\cdots,\lambda_d)$, alors il existe une
matrice
$M \in \on{GL}(d,\bhol{,F})$ telle que $(d_i)=M\Lambda (e_i)$ o\`u $(d_i)$ est une base de
$\dcris(V)$ et $(e_i)$ est une base de $\bhol{,F} \otimes_{\bplus_F} N_V$. 
Si $G_n$ est la matrice de l'action de $\gamma_n$ dans
la base vecteur colonne $(e_i)$: $\gamma_n(e_i)=G_n(e_i)$,
alors $\gamma_n(M\Lambda)G_n=M\Lambda$ ou encore
\[G_n=\gamma_n(\Lambda^{-1}M^{-1})M\Lambda
=\Lambda^{-1}D_n\gamma_n(M^{-1})M\Lambda 
\quad \text{et donc}\quad
\Lambda G_n \Lambda^{-1} = D_n \gamma_n(M^{-1})M 
\] o\`u l'on a 
pos\'e $D_n=\Lambda\gamma_n(\Lambda^{-1})$.
Remarquons que, comme $M \in \on{GL}(d,\bhol{,F})$, on a $\gamma_n(M^{-1})M=\on{Id}
\mod{\phi^{n-1}(q)}$. De plus, $D_n$ est une matrice diagonale et comme 
$\lambda_j=\pi^{\beta_{0,j}}(q/p)^{\beta_{1,j}}\cdots=(\phi^{n-1}(q))^{\beta_{n,j}}r_{n,j}$, 
o\`u $r_{n,j}$ est inversible modulo $\phi^{n-1}(q)$, on voit que
$\lambda_j/\gamma_n(\lambda_j)=\chi(\gamma_n)^{-\beta_{n,j}} \mod{\phi^{n-1}(q)}$. 

Les coefficients du polyn\^ome caract\'eristique de $G_n$  sont ceux de $\Lambda
G_n \Lambda^{-1}$ et sont donc \'egaux modulo $\phi^{n-1}(q)$ \`a ceux de $D_n$. Pour montrer que
$\Gamma_n$ agit par $(D_n \mod{\phi^{n-1}(q)})$ sur une base convenable de
$V_n=N_V/\phi^{n-1}(q)=(\bhol{,F} \otimes_{\bplus_F} N_V)/\phi^{n-1}(q)$, il suffit donc de
montrer que l'action de $\Gamma_n$ est diagonalisable, puisque les valeurs propres 
(avec multiplicit\'es) de $G_n$ seront alors celles de $(D_n \mod{\phi^{n-1}(q)})$.

L'application $\iota_n = \theta\circ\phi^{-n}$ r\'ealise une injection de $V_n$
dans $(\Cp \otimes_{\Qp} V)^{H_F} \simeq \oplus_{j=1}^d 
\hat{F}_{\infty}(-r_j)$ ($\iota_n$ est injective car un
\'el\'ement de $N_V$ divisible par $\phi^{n-1}(q)$ dans $\bplus\otimes_{\Qp} V$ est
divisible par $\phi^{n-1}(q)$ dans $\dfont^+(V)$ et donc dans $N_V$). Notamment, l'op\'erateur 
$\nabla_V(\nabla_V+1)\cdots(\nabla_V+r)$ est nul sur $V_n$, et ce dernier
espace est donc somme directe des $V_n^{\nabla_V=-j}$, qui sont stables par
$\Gamma_n$, et sur lesquels $\Gamma_n$ agit via le produit de $\chi^{-j}$ par une
repr\'esentation de $\Gamma_n$ qui est triviale sur un sous-groupe ouvert. 
L'action de $\Gamma_n$ est donc diagonalisable (l'action d'un groupe 
ab\'elien fini l'\'etant toujours).
\end{proof}

\begin{prop}
L'action de $\Gamma_n$ sur $V_n=N_V/\phi^{n-1}(q)$ 
est diagonale dans une base convenable de ce
$F_n$-espace vectoriel et $\gamma \in \Gamma_n$ agit dans cette base par
$\on{Diag}(\chi(\gamma)^{-r_1},\cdots,\chi(\gamma)^{-r_d})$
o\`u $r_1 \leq \cdots \leq r_d$ sont les oppos\'es des poids de Hodge-Tate de $V$.
\end{prop}

\begin{proof}
La proposition pr\'ec\'edente montre que dans une certaine base
de $V_n$, l'action de $\gamma \in \Gamma_n$ est donn\'ee par
$\on{Diag}(\chi(\gamma)^{-s_1},\cdots,\chi(\gamma)^{-s_d})$
o\`u les $s_1 \leq \cdots \leq s_d$ sont des entiers (ce sont les $\beta_{n,j}$). 
De plus, $\iota_n$ r\'ealise une
injection de $V_n$ dans
$(\Cp \otimes_{\Qp} V)^{H_F}_{\mathrm{fini}}=\dsen(V)$, qui est un $F_{\infty}$-espace vectoriel
muni d'une action r\'esiduelle de $\Gamma_F$ telle que $\nabla_V$ agisse, dans une bonne base,
par $\on{Diag}(-r_1,\cdots,-r_d)$

Dans cette base, l'action de $\gamma \in \Gamma_F$ est donn\'ee par une matrice \`a
coefficients dans $F_k$ pour $k \gg 0$, et le $F_k$-espace vectoriel engendr\'e par
cette base, $\dsen^k(V)$ est stable par $\Gamma_F$ (et donc $\Gamma_n$)
et v\'erifie $F_{\infty}\otimes_{F_k}
\dsen^k(V)=\dsen(V)$. De plus, quitte \`a augmenter $k$, l'image de $V_n$ dans 
$\dsen(V)$ est contenue dans $\dsen^k(V)$.

En exponentiant la connection $\nabla_V$, on voit que
quitte \`a augmenter $k$, l'action de
$\gamma \in \Gamma_k$ sur $\dsen^k(V)$ est diagonale et donn\'ee par
$\on{Diag}(\chi(\gamma)^{-r_1},\cdots,\chi(\gamma)^{-r_d})$.
La proposition r\'esulte alors du lemme
suivant que nous d\'emontrons ensuite:
\begin{lemm}
Si $V_n$ est un $F_n$-espace vectoriel muni d'une action de $\Gamma_n$, 
telle que dans une certaine base, $\gamma \in \Gamma_n$ agit par
$\on{Diag}(\chi(\gamma)^{-s_1},\cdots,\chi(\gamma)^{-s_d})$,
avec $s_1 \leq \cdots \leq s_d$, 
qui s'injecte de mani\`ere $\Gamma_n$-\'equivariante dans
$V_k$, un $F_k$-espace vectoriel muni d'une action de $\Gamma_n$ 
telle que l'action de $\gamma \in \Gamma_k$ est 
donn\'ee dans une certaine base par 
$\on{Diag}(\chi(\gamma)^{-r_1},\cdots,\chi(\gamma)^{-r_d})$
avec $r_1 \leq \cdots \leq r_d$, alors $r_i=s_i$.
\end{lemm}

En regardant les espaces propres pour $\nabla_V$ et en tordant, on se ram\`ene \`a montrer
que si $V_n$ un $F_n$-espace vectoriel muni d'une action triviale de $\Gamma_n$, qui
s'injecte de mani\`ere $\Gamma_n$-\'equivariante
dans un $F_k$-espace vectoriel $V_k$ muni d'une action de $\Gamma_n$,
alors $\dim_{F_n}(V_n) \leq \dim_{F_k}(V_k)$, ce qui suit du fait 
(que nous allons d\'emontrer) que
l'application naturelle $f: F_k \otimes_{F_n} V_n \ra V_k$ est injective. 
D\'emontrons cela: si ce n'\'etait pas le
cas, il existerait une relation minimale $\sum \mu_j f(v_j)=0$, et comme
$\Gamma_n$ agit trivialement sur $V_n$, on voit que $\Gamma_n$ fixe $\mu_j/\mu_1$ 
(comme la relation est minimale) et donc que
$\mu_j/\mu_1 \in F_n$ ce qui fait qu'un \'el\'ement 
du noyau de $f$ est d\'ej\`a nul dans
$F_k \otimes_{F_n} V_n$.

Ceci ach\`eve de d\'emontrer la proposition.
\end{proof}

\subsection{L'inclusion $\dfont^+(V) \subset \bplus \otimes_{\Qp} V$}
Le corollaire suivant est une cons\'equence des calculs du paragraphe \ref{gammaqn}
(pour $n \geq 1$; pour $n=0$, c'est une cons\'equence du fait que $\Gamma_F$ agit trivialement 
sur $N_V /\pi N_V$):
\begin{coro}
Dans les notations de \ref{invrigcris}, on a
$\beta_{0,j}=0$ et $\beta_{n,j}=r_j$ pour tout $n \geq 1$.
\end{coro}

On utilise ce corollaire pour montrer le r\'esultat principal de ce chapitre:
\begin{theo}\label{boundcomp}
Si $T$ est une $\Zp$-repr\'esentation de $G_F$, telle que $V = \Qp \otimes_{\Zp} T$ 
est cristalline \`a poids de Hodge-Tate entre $-r$ et $0$,
et $N_T$ est le $\aplus_F$-module dont on a rappel\'e la d\'efinition plus haut,
alors $\pi^r \aplus \otimes_{\Zp} T \subset \aplus \otimes_{\aplus_F} N_T$. En
particulier,
$\pi^r \aplus \otimes_{\Zp} T \subset \aplus \otimes_{\aplus_F} \dfont^+(T)$.
\end{theo}

\begin{proof}
On voit qu'il suffit de montrer que 
$\pi^r \bplus \otimes_{\Qp} V \subset \bplus \otimes_{\bplus_F}
N_V$. Montrons tout d'abord que $q^r$ tue $N_V/\phi^*N_V$. Comme $\phi^*(\bhol{,F}\otimes_F
\dcris(V))=\bhol{,F}\otimes_F
\dcris(V)$, on voit que $N_V / \phi^* N_V$ est tu\'e par 
$\prod_{i=1}^d \lambda_i/\phi(\lambda_i)$ qui est
une puissance de $q$. D'autre
part, si $x\in N_V$, alors $t^r x \in \bhol{,F}\otimes_F
\dcris(V) \subset \phi^* \bhol{,F} \otimes_{\bplus_F} N_V$, ce qui fait que $t^r$ tue 
$\bhol{,F} \otimes_{\bplus_F} N_V / \phi^*(\bhol{,F} \otimes_{\bplus_F} N_V)=
N_V/\phi^*N_V$ et donc $q^r$ tue $N_V/\phi^*N_V$.

Ceci montre que $q^r$ tue $N_V/\phi^*N_V$, ce qui revient \`a dire que si $P$ est la matrice de
$\phi$ dans une base de $N_V$, alors $q^r P^{-1} \in \on{M}(d,\bplus_F)$. Si $M$ est la matrice de
cette base dans une base de $V$, alors $\phi(M)=PM$ et $\phi(\pi^r M^{-1})=(\pi^r M^{-1})(q^r
P^{-1})$ ce qui fait que, par l'argument de r\'egularisation de \cite{Moi}, on a $\pi^r M^{-1} \in
\on{M}(d,\bplus)$ et donc $\pi^r \bplus \otimes_{\Qp} V \subset \bplus \otimes_{\bplus_F} N_V$.
\end{proof}

\begin{rema}
Le r\'esultat pr\'ec\'edent est proche d'\^etre optimal. On peut par exemple montrer que
$\bplus \otimes_{\Qp} V = \bplus \otimes_{\bplus_F} \dfont^+(V)$ si et seulement si la restriction de
$V$ \`a $H_F$ est non-ramifi\'ee.
\end{rema}

\section{Passage \`a la limite}
\subsection{L'image de $\dfont^+(T)$ dans $(\aplus/p^n \otimes_{\Zp} T)^{H_F}$}
Le foncteur $T \mapsto \dfont^+(T)$ n'est pas exact, mais en utilisant les r\'esultats de la
section pr\'ec\'edente, on peut contr\^oler l'image de $\dfont^+(T)$ dans $(\aplus/p^n
\otimes_{\Zp} T)^{H_F}$:
\begin{prop}\label{image}
Si $T$ est une $\Zp$-repr\'esentation de hauteur finie de $G_F$, et $\lambda \in \aplus_F$
est tel que $\lambda \aplus \otimes_{\Zp} T \subset \aplus \otimes_{\aplus_F} \dfont^+(T)$, 
et si $\mu \in \aplus_F$
est tel que $\mu \aplus \otimes_{\Zp} T \subset \aplus \otimes_{\aplus_F} N_T$, 
alors
\begin{enumerate}
\item l'image de $\dfont^+(T)$ dans $(\aplus/p^n \otimes_{\Zp} T)^{H_F}$ contient 
$\lambda (\aplus/p^n \otimes_{\Zp} T)^{H_F}$;
\item l'image de $N_T$ dans $(\aplus/p^n \otimes_{\Zp} T)^{H_F}$ contient 
$\mu (\aplus/p^n \otimes_{\Zp} T)^{H_F}$;
\item $\dfont^+(T) \simeq \projlim_n (\aplus/p^n \otimes_{\Zp} T)^{H_F}$;
\item $\dfont^+(T) \ra (\aplus/p^n \otimes_{\Zp} T)^{H_F}$ est surjective si $n\gg 0$.
\end{enumerate}
\end{prop}

\begin{proof}
Remarquons que $\lambda\in\aplus_F$ tel que $\lambda \aplus \otimes_{\Zp} T \subset 
\aplus \otimes_{\aplus_F} \dfont^+(T)$ existe toujours. Si $V$ est cristalline
\`a poids de Hodge-Tate dans $[-r;0]$, on vient de voir que l'on peut prendre $\lambda=\mu=\pi^r$.

Commen\c{c}ons par montrer le (1).
Soit $\beta\in(\aplus/p^n \otimes_{\Zp} T)^{H_F}$, que l'on rel\`eve en $\hat{\beta} \in
\aplus\otimes_{\Zp} T$. On a pour tout $h\in H_F$, $h(\hat{\beta})-\hat{\beta} \in p^n (\aplus
\otimes_{\Zp} T)$. On a aussi $\lambda \hat{\beta} \in \aplus \otimes_{\aplus_F} \dfont^+(T)$, et donc
\[ h(\lambda \hat{\beta})-\lambda \hat{\beta} \in p^n (\aplus
\otimes_{\Zp} T) \cap  \aplus \otimes_{\aplus_F} \dfont^+(T) \subset p^n \aplus \otimes_{\aplus_F} 
\dfont^+(T) \]
ce qui fait que si l'on \'ecrit $\lambda \hat{\beta} = \sum y_i \otimes d_i \in
\aplus \otimes_{\aplus_F} \dfont^+(T)$, alors $h(y_i)-y_i \in p^n \aplus$. Comme $\aplus_F \ra
(\aplus/p^n)^{H_F}$  est surjective, c'est que $y_i \in \aplus_F+p^n \aplus$, et donc que
$\lambda \hat{\beta} \in \dfont^+(T) + p^n \aplus \otimes_{\aplus_F} \dfont^+(T)$. Cela montre le (1).

Le (2) se d\'emontre exactement de la m\^eme mani\`ere, en remapla\c{c}ant $\dfont^+(T)$ par $N_T$.

Montrons le (3) de la proposition: soit $\{\beta_n\}$ une suite de $\projlim_n (\aplus/p^n
\otimes_{\Zp} T)^{H_F}$. Il existe une suite $\alpha_n$ de $\dfont^+(T)$ telle que l'image de
$\alpha_n$ est $\lambda\beta_n$. Cela implique imm\'ediatement que la suite $\alpha_n$
converge pour la topologie $p$-adique vers $\alpha \in \dfont^+(T)$ tel que l'image de
$\alpha$ est $\lambda\beta_n$. On peut \'ecrire, dans $\aplus\otimes_{\Zp} T$, $\alpha=\lambda
x_n + p^n y_n$ (o\`u $x_n \in \aplus \otimes_{\Zp} T$ rel\`eve $\beta_n$),
et la suite $x_n$ converge vers un \'el\'ement $x$ de $\aplus\otimes_{\Zp} T$
tel que $\alpha=\lambda x$. Ceci fait que $\alpha\in \dfont^+(T)$ a ses coordonn\'ees
divisibles par $\lambda$ et donc que $\alpha$ est divisible par $\lambda$ dans
$\dfont^+(T)$; $\alpha\lambda^{-1}$ s'envoie alors sur $\{\beta_n\}$, ce qui montre
le (3).

Passons maintenant au (4). Soit $D_n$ l'image de $(\aplus/p^n \otimes_{\Zp} T)^{H_F}$ dans
$(\eplus \otimes_{\Zp} T)^{H_F}$ et $D_{\infty}$ celle de $\dfont^+(T)$. On a clairement
$D_{\infty} \subset \cdots \subset D_n \subset \cdots \subset D_1$. Ils sont tous libres
de rang $d$ sur $\e_F^+$, dont les id\'eaux sont les $\pi^i$ avec $i\in\mathbf{N}$, et 
un petit argument de d\'eterminants montre que la suite d\'ecroissante $D_n$
doit stationner (puisqu'elle est born\'ee inf\'erieurement par $D_{\infty}$), disons que
$D_{n+1}=D_n$ pour $n\geq n_0$. Cela veut dire que $(\aplus/p^{n+1} \otimes_{\Zp} T)^{H_F} \ra
(\aplus/p^n \otimes_{\Zp} T)^{H_F}$ est surjective pour $n\geq n_0$, puisque c'est vrai modulo
$p$.  Soit $x\in (\aplus/p^n
\otimes_{\Zp} T)^{H_F}$, on peut donc le relever en un \'el\'ement de $\projlim_n (\aplus/p^n
\otimes_{\Zp} T)^{H_F}$, et par le (3), en $\hat{x}\in\dfont^+(T)$.
\end{proof}

\subsection{Construction de $N_T$}\label{col1}
L'objet de cette section est de d\'emontrer la conjecture de Fontaine sur les 
limites de repr\'esentations
cristallines de $G_F$, c'est-\`a-dire le th\'eor\`eme suivant:
\begin{theo}\label{cristheo}
Si $T$ est une $\Zp$-repr\'esentation de $G_F$, telle qu'il existe $b \geq a \in \NN$ et une suite
$\{T_i\}_i$ de $\Zp$-repr\'esentations cristallines de $G_F$ \`a poids de Hodge-Tate dans $[a;b]$
v\'erifiant $T/p^i \simeq T_i/p^i$, alors $T$ est cristalline.
\end{theo}

Tout d'abord, en tordant $T$ et les $T_i$, on se ram\`ene au cas o\`u les $T_i$ ont leurs poids de 
Hodge-Tate dans $[-r;0]$ pour un $r \geq 0$. Dans tout ce paragraphe, $N_{T_i}$ d\'enotera le module
de Wach associ\'e \`a $T_i$ dont on a rappel\'e la construction plus haut, et qui fait l'objet du
th\'eor\`eme \ref{boundcomp}. En particulier, on dispose du lemme suivant qui rassemble l'information
dont on aura besoin:

\begin{lemm}
Le $\aplus_F$-module $N_{T_i}$ est un sous $\aplus_F$-module de $\dfont^+(T_i)$, stable par $\phi$
et $\Gamma_F$, tel que pour tout $\gamma \in \Gamma_F$, on ait $(\gamma-1)N_{T_i} \subset
\pi N_{T_i}$. De plus $\pi^r \aplus \otimes_{\Zp} T_i
\subset \aplus \otimes_{\aplus_F} N_{T_i}$ et
l'image de $N_{T_i}$ dans $(\aplus/p^n \otimes_{\Zp} T_i)^{H_F}$ contient 
$\pi^r (\aplus/p^n \otimes_{\Zp} T_i)^{H_F}$. 
\end{lemm}

Remarquons que $N_{T_n}$ s'envoie naturellement dans 
$(\aplus/p^n \otimes_{\Zp} T_n)^{H_F} \simeq (\aplus/p^n \otimes_{\Zp} T)^{H_F}$, et
donc dans $(\aplus/p^m \otimes_{\Zp} T)^{H_F}$ si $n \geq m$. On va d'abord montrer la
proposition suivante, qui permettra de recoller les $N_{T_i}$ en un $N_T$.

\begin{prop}
Si $T_1$ et $T_2$ sont deux r\'eseaux de deux repr\'esentations cristallines \`a poids de
Hodge-Tate dans $[-r;0]$, tels que $T_1 / p^n = T_2 / p^n$
avec $n \in \NN$ v\'erifiant $p^{n-1}(p-1) \geq r+1$, alors
$N_{T_1}$ et $N_{T_2}$ ont m\^eme image dans $(\aplus/p^n \otimes_{\Zp/p^n} T_i/p^n)^{H_F}$.
\end{prop}

\begin{proof}
L'id\'ee de la d\'emonstration est que les images de $N_{T_1}$ et $N_{T_2}$ dans
$D_n = (\aplus/p^n \otimes_{\Zp/p^n} T_i/p^n)^{H_F}$ sont deux $\aplus_F/p^n$-modules libres,
stables par $\Gamma_F$, tels que $\Gamma_F$ agit trivialement sur $N_{T_1}/\pi$ et
$N_{T_2}/\pi$, et tels que pour tout $x \in D_n$, il existe $m \leq r$ tel que 
$\pi^m x \in N_{T_i}/p^n$ (par le lemme pr\'ec\'edent). 
Si $n=+\infty$, alors on a vu en \ref{constnw} qu'un tel objet est unique.

Dans notre cas, soit $x \in N_{T_1}/p^n \setminus \pi  N_{T_1}/p^n$. Il existe 
$m \leq r$ tel que $\pi^m x \in N_{T_2}/p^n \setminus \pi  N_{T_2}/p^n$, 
et en regardant l'action de $\Gamma_F$, on voit que
l'on doit avoir $\gamma(\pi^m)/\pi^m = 1 \mod{(p^n,\pi)}$, et donc
$\chi(\gamma)^m = 1 \mod{p^n}$ pour tout $\gamma \in \Gamma_F$. Cela implique que
$p^{n-1}(p-1)$ divise $m \leq r$, mais on a suppos\'e que $p^{n-1}(p-1) \geq r+1$. C'est
donc que $m=0$, et que $N_{T_1}/p^n \subset N_{T_2}/p^n$. Par sym\'etrie, on en conclut
que $N_{T_1}/p^n \subset N_{T_2}/p^n$.
\end{proof}

\begin{rema}
La borne $p^{n-1}(p-1) \geq r+1$ est essentielle, on peut facilement construire des
contre-exemples \`a la proposition pr\'ec\'edente si $p^{n-1}(p-1) \leq r$.
\end{rema}

\begin{prop}\label{colmez}
Si $N(T)= \projlim_m N_{T_m}/p^m$, alors $N(T)$ est un sous-$\aplus_F$-module 
libre de rang $d$ de $\dfont(T)$, stable par $\phi$ et $\Gamma_F$ et tel que $\Gamma_F$
agisse trivialement sur $N(T)/ \pi N(T)$. 
\end{prop}

La proposition ci-dessus implique le th\'eor\`eme \ref{cristheo}, puisque $N(T)$  est un sous
$\aplus_F$-module de $\dfont(T)$, libre de rang $d$, stable par $\phi$ et $\Gamma_F$, et tel que
$\Gamma_F$ agisse trivialement sur $N(T) / \pi N(T)$. En particulier, $T$ est de hauteur
finie et cristalline (et bien s\^ur $N(T)=N_T$). Reste donc \`a montrer la proposition
\ref{colmez}.

\begin{proof}
Il est clair que $N(T)$ s'identifie \`a un sous-$\aplus_F$-module de
$\projlim_m (\mathbf{A} /p^m \otimes_{\Zp} T)^{H_F} \simeq \dfont(T)$, 
stable par $\phi$ et $\Gamma_F$ et tel que $\Gamma_F$ agisse 
trivialement sur $N(T)/ \pi N(T)$. 
Reste \`a montrer que $N(T)$ est un $\aplus_F$-module libre de
rang $d$. Par la proposition pr\'ec\'edente, d\`es que $p^{n-1}(p-1) \geq r+1$, 
$N(T)/p^n \simeq N_{T_n}/p^n$ est un 
$\aplus_F/p^n$-module libre de rang $d$. 

L'application naturelle de $N(T)$ dans $N(T)/p^n$
est surjective et $N(T)/p^n$ est libre de rang $d$, ce qui
fait qu'une base de $N(T)/p^n$ se rel\`eve, par le lemme de Nakayama ($N(T)$ est s\'epar\'e
complet pour la topologie $p$-adique, car chaque $N_{T_n}$ l'est), en une famille libre et
g\'en\'eratrice de $N(T)$, qui est donc libre de rang $d$.
\end{proof}

Ceci ach\`eve de d\'emontrer le th\'eor\`eme \ref{cristheo}. 
\begin{coro}
Si $T$ est une $\Zp$-repr\'esentation de $G_F$, telle qu'il existe $r \in \NN$ et une suite
$\{T_i\}_i$ de $\Zp$-repr\'esentations cristallines de $G_F$ \`a poids de Hodge-Tate dans $[a_i;b_i]$
avec $b_i-a_i \leq r$,
v\'erifiant $T/p^i \simeq T_i/p^i$, alors il existe un caract\`ere $\eta:\Gamma_F \ra \Zp^*$ tel
que $T(\eta^{-1})$ est cristalline.
\end{coro}

\begin{rema}
Si l'on regarde la d\'emonstration du th\'eor\`eme \ref{cristheo}, on voit que l'on a en fait
d\'emontr\'e que si $T$ est une $\Zp$-repr\'esentation de $G_F$, telle qu'il existe $a \in \NN$ et une
suite $\{T_i\}_i$ de $\Zp$-repr\'esentations cristallines de $G_F$, \`a poids de Hodge-Tate dans
$[a;a+p^{i-2}(p-1)-1]$ pour $i \gg 0$,
v\'erifiant $T/p^i \simeq T_i/p^i$, alors $T$ est cristalline.

En effet, on voit dans ce cas que $N_{T_i} / p^i = N_{T_{i+1}} / p^i$ et donc que l'on peut
recoller les $N_{T_i}$ comme ci-dessus.
\end{rema}

\section{L'identification $N_V/\pi N_V \simeq \dcris(V)$}
\subsection{Construction de $\phi$-modules filtr\'es}
Le but de ce chapitre est de donner une d\'emon\-stration du th\'eor\`eme suivant:
\begin{theo}\label{crisnv}
Si $V$ est une repr\'esentation cristalline positive de $G_F$, 
et si $N_V$ est le 
$\bplus_F$-module dont on rappel\'e la construction plus haut,  
muni de la filtration $\on{Fil}^i N_V = \{ x \in N_V, \phi(x) \in
q^i N_V \}$,
alors $\dcris(V) \subset \bhol{,F} \otimes_{\bplus_F} N_V$ et
l'application naturelle $\dcris(V) \ra N_V / \pi N_V$ est 
un isomorphisme de $\phi$-modules filtr\'es.
\end{theo}

\begin{proof}
Pour montrer que cette application est bijective, il suffit de voir que $\pi N_V \cap
\dcris(V) =\{0\}$, ce qui suit du fait que si $x \in \dcris(V)$, alors $\gamma(x)=x$ tandis
que si $x \in \pi N_V$, alors  $\gamma(x)-\chi(\gamma)x \in \pi^2 N_V$, et donc si 
$x \in \pi N_V \cap \dcris(V)$, alors $x\in\pi^2
N_V$. En it\'erant ce proc\'ed\'e, on voit que $x \in \cap_{j=0}^{+\infty } \pi^j N_V =
\{0\}$.

Il est clair que cette application est un isomorphisme de $\phi$-modules. Enfin, 
un \'el\'ement de $\btrigplus{}$ est divisible par $\pi/\pi_1$ dans
$\bmax^+$ si et seulement s'il l'est dans $\btrigplus{}$, ce qui fait que
$x\in \on{Fil}^i(\btrigplus{})$ si et seulement si $\phi(x) \in q^i \btrigplus{}$. 
On en d\'eduit que si $x \in \on{Fil}^i(\btrigplus{} \otimes_{\Qp} V)^{G_F}$, 
alors son image dans $\bhol{,F} \otimes_{\bplus_F} N_V$ v\'erifie 
\[ \phi(x) \in (q^i \btrigplus{} \otimes_{\Qp} V) \cap (\bhol{,F} \otimes_{\bplus_F} N_V) =
q^i \bhol{,F} \otimes_{\bplus_F} N_V. \]
En effet, $\bplus[1/\pi] \otimes_{\Qp} V = \bplus[1/\pi] \otimes_{\bplus_F} N_V$ et donc
$\btrigplus{}[1/\pi] \otimes_{\Qp} V = \btrigplus{}[1/\pi] \otimes_{\bplus_F} N_V$, ce qui fait que
$q^i (\btrigplus{} \otimes_{\Qp} V) \cap (\bhol{,F} \otimes_{\bplus_F} N_V) 
= q^i \bhol{,F} \otimes_{\bplus_F} N_V$.

Ceci montre que l'image de $\on{Fil}^i \dcris(V)$ est incluse dans $\on{Fil}^i
(N_V / \pi N_V)$. R\'eciproquement, si $y \in \dcris(V) 
\subset \bhol{,F} \otimes_{\bplus_F} N_V$ a la
propri\'et\'e que
$\phi(y) \in q^i \bhol{,F} \otimes_{\bplus_F} N_V$, alors 
$\phi(y) \in q^i (\btrigplus{} \otimes_{\Qp} V)$ et
donc $y \in \on{Fil}^i \dcris(V)$.
\end{proof}

\subsection{Construction de $\dfont^+(T)$ dans le cas $r \leq p-2$}\label{col2}
Si $r_d-r_1 \leq p-2$, on peut donner une construction explicite de $\dfont^+(T)$, 
qui montre directement que $\aplus[1/\pi] \otimes_{\Zp} T = \aplus[1/\pi] \otimes_{\aplus_F} 
\dfont^+(T)$.
Ces r\'esultats 
(construction de $\Zp$-repr\'esentations cristallines associ\'ees \`a des r\'eseaux fortement
divisibles) sont d\'ej\`a connus \cite{L80,FL82,Wa97}, 
mais notre d\'emonstration est plus directe. Le lemme ci-dessous est imm\'ediat:

\begin{lemm}\label{procheqr}
Si $\mu=(p/(q-\pi^{p-1})) \in \aplus_F$, alors $\mu(0)=1$ 
($\mu$ est donc inversible dans $\aplus_F$), et 
$\mu^s q^s = p^s \mod{\pi^{p-1}}$.
\end{lemm}

\begin{theo}
Si $V$ est une repr\'esentation cristalline positive
de $G_F$, telle que la longueur de la filtration
sur $\dcris(V)$ est \'egale \`a $r \leq p-2$, et si $M$ est un r\'eseau 
fortement divisible de $\dcris(V)$, alors il existe un r\'eseau $T$ de $V$, et
une base de $\dfont(T)$ qui engendre un
$\aplus_F$-module $N_T$, tel que: $q^r N_T \subset \phi^*(N_T) \subset N_T$, l'action de
$\Gamma_F$ est triviale sur $N_T/\pi N_T$, et si l'on munit $N_V=\Qp\otimes_{\Zp}N_T$ de la
filtration d\'efinie par $\on{Fil}^i N_V = \{ x \in N_V, \phi(x) \in q^i N_V \}$, alors
$N_V/\pi N_V$ est isomorphe \`a $\Qp\otimes_{\Zp}M$ en tant que $\phi$-modules filtr\'es.
\end{theo}

\begin{proof}
Quitte \`a tordre $V$, on se ram\`ene au cas o\`u $r_1=0$ et donc $r_d \leq p-2$.
Si $M$ est un r\'eseau fortement divisible de $\dcris(V)$, alors dans une base adapt\'ee \`a
la filtration, la matrice de $\phi$ est \'egale \`a $P_0 A$ o\`u
$P_0=\on{Diag}(p^{r_1},\cdots,p^{r_d})$ et $A \in \on{GL}(d,\OO_F)$. On pose
$P=\on{Diag}(q^{r_1}\mu^{r_1},\cdots,q^{r_d}\mu^{r_d})A$ ce qui fait que $P$ est \'egale
modulo $\pi$ \`a la matrice de $\phi$ sur $M$. De plus, le lemme 
\ref{procheqr} montre que si $\gamma \in \Gamma_F$, alors $\gamma(P^{-1})P
\in \on{Id}+\pi^{p-1} \on{M}(d,\aplus_F)$.  On \'ecrira $\gamma(P^{-1})P
= \on{Id} + \pi^{p-1} Q$.

En particulier,
si $\gamma \in \Gamma_F$, alors l'application $f_{\gamma}$ 
d\'efinie par $f_{\gamma}(H)=\gamma(P^{-1})\phi(H)P$ envoie
$\on{Id}+\pi^{p-1}\on{M}(d,\aplus_F)$ dans 
$\on{Id}+\pi^{p-1}\on{M}(d,\aplus_F)$ et une deuxi\`eme application du lemme
\ref{procheqr} montre que $f_{\gamma}$ 
est contractante pour la topologie $(p,\pi)$-adique, et y admet
donc un unique point fixe $G(\gamma)$, qui satisfait alors $\phi(G(\gamma))P=\gamma(P)G(\gamma)$.

On peut d\'efinir un
$(\phi,\Gamma)$-module $D$, en faisant agir $\phi$ et $\gamma \in \Gamma_F$ par $P$ et
$G(\gamma)$. Ce $(\phi,\Gamma)$-module \'etale $D$ correspond \`a un r\'eseau $T$ d'une 
repr\'esentation  de hauteur finie, et le th\'eor\`eme \ref{crisnv} montre que cette
repr\'esentation est $V$ (puisque le foncteur $\dcris(\cdot)$ est pleinement fid\`ele), car $N_V/\pi
N_V$ est un $\phi$-module filtr\'e sur lequel $\phi$ agit  par la m\^eme matrice que sur $M$, ce qui
montre que $N_V/\pi N_V \simeq \Qp \otimes_{\Zp}M$  en tant que
$\phi$-modules filtr\'es (la filtration sur $N_V$ \'etant d\'efinie via l'action de $\phi$).
\end{proof}

\backmatter

\end{document}